\documentclass[4pt]{article}
\usepackage{amsmath,amsfonts,amsthm,amssymb,amscd}
    \usepackage[cal=boondoxo]{mathalfa} 
    \usepackage{hyperref} 

\newtheorem{thm}{Theorem}[section]

\newtheorem*{thm*}{Theorem}
\newtheorem*{corr*}{Corollary}
\newtheorem{lemma}[thm]{Lemma}
\newtheorem{prop}[thm]{Proposition}
\newtheorem*{prop*}{Proposition}

\newtheorem{corr}[thm]{Corollary}
\theoremstyle{definition}
\newtheorem{dfn}[thm]{Definition}

\newtheorem{exmples}[thm]{Examples}

\theoremstyle{remark}
 
\newtheorem*{rmq}{\textit{Remark}}
\newtheorem{rmk}[thm]{\textit{Remark}}
\newtheorem{rmks}[thm]{\textit{Remarks}}
\renewcommand{\proof}{\noindent\textit{Proof}\/: \,\,}

\def\C{{\mathbf{C}}}
\newcommand{\Q}{{\mathbf{Q}}}
\newcommand{\R}{{\mathbf{R}}}

\newcommand{\bH}{{\mathbf{H}}}

\newcommand{\bS}{{\mathbf{S}}}

\newcommand\germ[1]{{\mathfrak{#1}}}
\newcommand{\eh}{\germ h}%

\renewcommand{\AA}{{\mathcal A}}

\newcommand\FF{{\mathcal F}} 
\newcommand\GG{{\mathcal G}} 
\newcommand\HH{{\mathcal H}}

\newcommand\MM{{\mathcal M}} 
\newcommand\OO{{\mathcal O}} 
\newcommand\TT{{\mathcal T}}

\newcommand{\comp}{\raise1pt\hbox{{$\scriptscriptstyle\circ$}}}
 
\newcommand{\underint}{\int_{\raise-4pt\hbox{\hskip-8pt $-$}}}
\newcommand{\overint}{\int^{\raise3pt\hbox{\hskip-7pt $-$}
}\hskip -4pt}

\def\lset{\{}  
\def\rset{\}}  
\def\set#1{\lset#1\rset} 
\def\st{\mid}   
\def\sett#1#2{\lset #1 \st #2 \rset}

\newcommand\Tr{{}^{\mathsf{T}}\kern-0.9pt} 

\newcounter{lijstc}
{\end{list}}
{\end{list}}

\def\mapright#1{\mathop{\vbox{\ialign{
                                ##\crcr
    ${\scriptstyle\hfil\;\;#1\;\;\hfil}$\crcr
 \noalign{\kern2pt\nointerlineskip}
    \rightarrowfill\crcr}}\;}}

\def\mapleft#1{\mathop{\vbox{\ialign{
                                ##\crcr
    ${\scriptstyle\hfil\;\;#1\;\;\hfil}$\crcr
 \noalign{\kern2pt\nointerlineskip}
    \leftarrowfill\crcr}}\;}}

\newcommand\rarrow[3]{\smash{\mathop{\hbox to#3{\rightarrowfill}}\limits
^{\scriptstyle#1}_{\scriptstyle#2}}}
 
\newcommand\larrow[3]{\smash{\mathop{\hbox to#3{\leftarrowfill}}\limits
^{\scriptstyle#1}_{\scriptstyle#2}}}

\renewcommand\emptyset{\varnothing}
\renewcommand\setminus{-}
\def\into{\hookrightarrow}
\def\onto{\twoheadrightarrow}

\newcommand\su[1]{\operatorname{SU}({#1})}

\newcommand\SL[1]{\operatorname{SL}({#1})}
\newcommand\smpl[1]{\operatorname{Sp}({#1})}

\def\half{\frac{1}{2}}
\def\Res{\operatorname{Res}}
\def\vhs{\text{variation of Hodge structure}}
\def\gen{\text{\rm gen}}
\def\jac{\text{\rm Jac}}
\def\bmy{\beta}
\def\rank{\operatorname{rank}}

\title{Inequalities for semi-stable  surface fibrations, and their relation to the Coleman-Oort conjecture}
\author{Chris Peters}
\date{}

\begin{document}
\maketitle
\section*{Introduction}

\subsection*{Background} 
Arakelov showed \cite{arakelov} that over a fixed curve there are at most finitely many non-constant fibrations  in  curves of genus $\ge 2$ for which its singular fibers appear over a fixed set in the base.
 To prove  this he first showed that these fibrations appear in a finite number of families (boundedness) and secondly that these are all rigid.
Boundedness follows from a bound  involving  the relative canonical bundle. This bound is  the "classical" Arakelov inequality.

The moduli space $\MM_g$ of genus $g$ curves is quasi-projective and admits a natural compactification,
 $\overline{\MM_g}$, the Deligne-Mumford  compactification. See e.g. \cite[Chapter X]{curves}. 
It  is constructed by adding  points to the boundary
which  correspond  to (isomorphism classes of) \emph{stable curves}.  These are (in general reducible) non-reduced curves 
 having ordinary double points and   each smooth rational component 
 meets  the remaining components in at least $3$ points. Equivalently: the fibration should first of all be 
  relatively minimal, i.e. no fiber contains a \emph{$(-1)$-curve}\footnote{I.e., a smooth rational curve  of self-intersection $(-1)$}, and secondly,
  no fiber should contain an  ADE-configuration.\footnote{I.e. a tree of  smooth rational curves of self-intersection $-2$ resolving a rational double point.} 
  So, in order to have a natural morphism of $B$ into  $\overline{\MM_g}$ all fibers should be
 stable. Clearly, this is not always the case, but there is a procedure ("semi-stable reduction") explained   e.g. in \cite[III.10]{4authors},  which produces
 a fibration over a finite (ramified) covering of $B$ whose fibers  in general not yet  stable but\emph{ semi-stable}: the fibers
 are reduced and have only double points and do not contain $(-1)$--curves but might contain $(-2)$--curves. The latter can be blown down in the surface which produces rational double points and the resulting family has stable fibers and hence its base maps to $\overline{\MM_g}$.
In view of these considerations, one may often reduce the study of curve fibrations
over curves to  semi-stable fibrations. I shall assume this for the remainder of this introduction.

Since   the Arakelov  bound is essentially Hodge theoretic in nature, 
 it is more natural to formulate it
in terms of variations of Hodge structure. The weight one case correspond to families of curves through their first
cohomology groups.  This link leads to an  an important observation: the cohomology of a  singular semi-stable  fiber  $C$ 
may have the
Hodge structure of a smooth fiber; this corresponds to  its Jacobian $\jac(C)$ being a polarized abelian variety. 
In our setting (of semi-stable fibrations) this happens precisely when $C$ is a tree of  of smooth curves, say $C_j$, $j=1,\dots,N$ and then
  $\jac(C)=\prod_{j=1}^N \jac(C_j)$.   Such a curve as well as its Jacobian is  called \emph{of compact type}.
 This shows that it is natural to separate the critical locus 
$\Delta$ of $f$ in two: the set of points $\Sigma_c$ over which the fiber is of compact type, and the remainder   $\Sigma$. Whereas   Arakelov's original 
 bound is expressed in terms of $\Delta$, the Hodge theoretic bound uses $\Sigma$. It involves  the relative canonical sheaf\footnote{See Section~\ref{sec:propdev} where  the definition and basic properties are given.} $\omega_{X/B}$ and can be expressed as follows.
  \begin{equation}
\label{eqn:BasicArakIneq}
\deg (f_*\omega_{X/B} )  \le  \half g \cdot  (2g(B)-2+\#\Sigma),\, g=\text{fiber genus},\, g(B)=\text{base genus}.
\end{equation}
In this guise it has been shown by Faltings. Hodge theoretic generalizations to higher weights can be found in \cite{Petersiv,JostZuo}.
It has become common practice to call this inequality also  an "Arakelov inequality".

\subsection*{Special curve fibrations over a curve}

 An interesting special case arises when \emph{the Arakelov bound is attained}, i.e. when equality holds in \eqref{eqn:BasicArakIneq}. 
 Rather surprisingly, one has \cite{shimura}:
\begin{thm*}[A]  For a  non-isotrivial   genus $g\ge 5 $ fibration $f:X\to B$ of semi-stable  curves
 the Arakelov bound    \eqref{eqn:BasicArakIneq} can never be attained.
\end{thm*}

 The aim of this note is first of all to  give a simplified and self contained proof of Theorem (A). 
 By using  refined Chern class inequalities for non-compact surfaces obtained in  \cite{kob,cy} instead of  those from 
 \cite{mi}, the
 proof given here avoids the elaborate covering tricks employed in \cite{shimura}. This not only simplifies their proof but also  makes the basic idea of their proof  more transparent.
 This is the content of    Sects.~\ref{sec:Surfaces}--\ref{sec:mainresult}.
 
 The main idea of the proof is as follows. For surface fibrations in semi-stable curves Moriwaki obtained an
 inequality   \eqref{eqn:slope} for the so called slope. 
 The   Chern class inequality   cited above  
is a refinement (=Theorem~\ref{kobaresult}) of the classical  Bogomolov-Miyaoka-Yau inequality\footnote{In what follows this will be abbreviated as \emph{BMY--inequality}.}.
 If the Arakelov bound is attained one  can rewrite this as an inequality for the slope  in the reverse direction.
 Note that this inequality depends on the boundary of the open surface. One has to choose this boundary carefully
 and then the combined inequalities are seen to be  only  possible when the fibering curves have genus at most $4$.

At this point I should mention that  ideas inspired by the cited work \cite{shimura,shimura2} of X. Lu and K. Zuo  have been used in the reverse direction   to obtain  more refined slope inequalities, a result
that should be of  direct interest to surface geometers. 
See \cite{ineqs1,ineqs2,ineqs3,ineqs4}.

\subsection*{Moduli aspects}

 A further aim of this note is to  elucidate the  \emph{moduli aspects}\footnote{A modern and self-contained  reference  for moduli spaces of curves  is  
\cite{curves}.}  of Theorem (A), 
especially its relation  with the so-called Coleman-Oort conjecture \cite{coleman,oort} about the Torelli locus. 
Let me  explain this very briefly. As   noticed before,
the Arakelov inequality \eqref{eqn:BasicArakIneq} admits  a purely Hodge theoretical proof which holds
for a polarized weight one variation  of Hodge structure\footnote{The proof of the resulting inequality \eqref{eqn:Arak0}
 is reproduced below in Section~\ref{sec:HiggsAndArak}.} over a curve, or, equivalently,  for a family of  principally polarized Abelian varieties of dimension $g$ over a curve.
 
 The corresponding moduli space $\AA_g$ contains  many obvious subvarieties coming from other  Abelian varieties 
 with extra structure, for instance those coming from products of lower dimensional Abelian varieties. 
Such subvarieties belong to the class of what nowadays is called the class of 
\emph{special subvarieties}, by definition  images in $\AA_g$ of  \emph{Shimura varieties}. See \S~\ref{sec:SShimura} for
a brief introduction to the subject of Shimura varieties.   Shimura varieties of dimension $1$
are   called \emph{Shimura curves}. Note that contrary to some practices,  I  do not
consider them as embedded in $\AA_g$, but as coming from a Shimura datum. In Section~\ref{sec:ShimCurves}    I
explain Satake's classification of Shimura curves:    the non-compact
Shimura curves are precisely the \emph{modular curves} parametrizing elliptic curves with a fixed  level structure and 
  the compact
ones come from the quotient of the upper half plane by  an  arithmetic subgroup of  the group of units
of certain quaternion algebras.

These curves come up naturally because of the  following basic result,  
essentially due to E. Viehweg  and K. Zuo (see \cite{carshcurves}):
\begin{prop*}[=Prop.~\ref{RigidSHiCurves}] If equality holds in the Hodge theoretic version  \eqref{eqn:Arak0}
of the Arakelov inequality, the punctured curve $B\setminus \Sigma$ is a rigidly embedded Shimura curve.
The converse also holds.
\end{prop*}
It is well known that $\AA_g$ is quasi projective and contains  many special subvarieties. For instance special points are dense. The  many types of special subvarieties of  $\AA_g$ all  have been classified \cite{addington}. It follows from this that the rigidly embedded special curves
form a relatively small list:
\begin{enumerate}
\item fibered self-products of the modular family of elliptic curves. The corresponding embedded curves are called \emph{of Satake type};
\item    fibered self-products
of modular families of Abelian varieties over compact  Shimura curves. The corresponding embedded curves are called
\emph{of Mumford type} after \cite{mumford}.
\end{enumerate}    

Another interesting sublocus of $\AA_g$  comes from the moduli space  $\MM_g$ of (smooth projective) genus $g$ curves: the Torelli theorem  states that the period map, which associates to a smooth genus $g$ curve  $C$ its Jacobian $\jac (C)$,  defines  an injective morphism
\[
p: \MM_g \to \AA_g.
\]
For $g\ge 2$ the image of $p$ is not closed in $\AA_g$. As is well known \cite[p. 74]{mum2} the limits come from 
curves of compact type encountered before: their (generalized) Jacobian is a product of lower dimensional Abelian varieties.
This non-closedness explains why, traditionally,  the \emph{Torelli locus} $\TT_g$ is defined to be the \emph{closure} inside  $\AA_g$
of the image of the period map.  

For $g\ge 4$ the   locus  $\TT_g$ is not special in the above sense. Indeed,  it is believed that  positive dimensional special subvarieties  $S\subset \AA_g$ cannot  generically be contained in the Torelli locus\footnote{ $S\subset \AA_g$ is  \emph{generically contained in $\TT_g$} if it is contained in $\TT_g$ but also meets $p(\MM_g)$
in a Zariski open subset of $S$.}  unless $g$ is "small''.
 This belief goes under the name of  the \emph{Coleman-Oort conjecture}, see  \cite{coleman,oort}. In \cite{Torloc} 
 background for this conjecture is given as well as examples why this fails for low $g$.  

An important conclusion of the previous discussion is that    Shimura curves    which are rigidly embedded in $\AA_g$  are precisely the
ones for which the corresponding Hodge bundle  attains the Arakelov bound. Hence,   Theorem (A)
clearly is related to the Coleman-Oort conjecture. Unfortunately, $\AA_g$ is not a fine moduli space and so one cannot assume
that  there is a (tautological) family of genus $g$ curves over an embedded special curve. This explains why further work needs to be done. In  \cite{shimura2} the authors show that one in fact has:
 \begin{thm*}[B]  Suppose $g\ge 12$. Then no  rigidly embedded  special curve (necessarily  of Satake  or of Mumford type)
  can   generically be contained in  $\TT_g$. 
\end{thm*} 
In Section~\ref{sec:Fin}   I shall briefly outline how this can be shown by suitably modifying the proof of Theorem (A) as given in this note. These modifications
are however not straightforward and require new ideas which I don't elaborate on.

Finally, I also want to  remark  that  similar ideas enabled  X. Lu and K. Zuo   \cite{shimura3}  to verify the Coleman-Oort conjecture for certain higher dimensional Shimura varieties.

 \medskip

\textbf{Acknowledgements}: I thank Kang Zuo for spotting a gap in a preliminary  version of my proof  of Theorem A
and his hints as how to fill  it.  Furthermore, I want to express thanks to  the anonymous referee(s)
for spotting several inaccuracies and unclear passages in this   modified version.

 \section{Background from surface theory} \label{sec:Surfaces}
 \subsection{Proportionality deviation for semi-stable fibrations} \label{sec:propdev}

 Let $X$ be a minimal compact complex algebraic surface and suppose that  
\begin{equation}\label{eqn:fibration}
f:X \to B,\quad \text{\rm genus } B= b,
\end{equation}
is a fibration in genus $g\ge 2$ semi-stable curves. Denote the fibre over $s$  by $X_s$; $\jac (X_s)$ denotes the generalized Jacobian of $X_s$. As usual, put 
\[
\omega_{X/B}= \OO_X(K_X\otimes f^*K_B^{-1}),\quad \text{the relative dualizing sheaf.}
\]
Recall  \cite[Theorem III, 18.2]{4authors} that a relatively minimal   fibration\footnote{I.e.,  there are no $(-1)$--curves in the fibres.} in genus $\ge 2$ curves $f:X\to B$ is isotrivial\footnote{I.e., locally constant over the non-critical locus of $f$.}  if and only if $\deg f_*\omega_{X/S}=0$. So, if $f$ is not isotrivial, one can introduce the  \emph{slope}
\begin{equation} \label{eqn:slope}
\lambda(f):= \omega^2_{X/B}/ \deg ( f_*\omega_{X/B}).
\end{equation}
Recall also that $f$ is called a \emph{Kodaira fibration} if  $f$ is smooth. In that case, the period map for $f$ sends $B$ to a compact curve in the moduli space $\MM_g$  of genus $g$ curves. For any  Kodaira fibration  the base genus  $g$
 has to be at least $2$. See \cite[V. 14]{4authors}. 

Relations between the various numerical invariants of a surface fibration are gathered in the next Lemma.

\begin{lemma} \label{basics}   For a relatively minimal  genus $g$ surface fibration\footnote{The fibers need not be semi-stable.}  $f:X\to B$, $g(B)=b$, one has\footnote{Standard notation is used here: $e$ denotes the topological Euler number and $c_1^2(X)=K^2_X$}:
\begin{enumerate}   \renewcommand{\labelenumi}{\rm (\roman{enumi})}
\item  $\deg ( f_*\omega_{X/B})=\frac 1 {12} (c_1^2(X) +e(X)) - (b-1)(g-1)$.
\item  $c_1^2(X)= \omega^2_{X/B} +8(g-1)(b-1).$
\item $
e(X)= 4(g-1)(b-1)+  \sum_{s\in B} \delta_s,
$
where
\begin{equation}\label{eqn:EgenEspec}
\delta_s := e(X_s)- e(X_\gen)\ge 0,\quad e(X_\gen) =2(1-g).
\end{equation}
\item $12 \deg f_*\omega_{X/B}= \omega_{X/B}^2 +\sum_s \delta_s.$
\item  Suppose, moreover, that $f$ is non-isotrivial and $g\ge 2$. Then 
 \begin{equation}\label{eqn:xiao}
 4 - \frac 4 g\le \lambda(f)\le 12.
\end{equation}
The leftmost inequality is Xiao's  \emph{slope inequality}. The right hand inequality becomes an equality,  $\lambda(f)=12$ if and only if $f$ is a Kodaira fibration. 
\end{enumerate}
\end{lemma}
\proof (i) is a direct consequence  of Riemann-Roch and the Leray spectral sequence.\\
(ii) follows from the definition of $\omega_{X/B}$.\\
(iii) See   \cite[Prop. III, 11.4]{4authors}. \\
(iv) Follows directly from (i),(ii), (iii).\\
(v) See  \cite{xiao}. The remaining assertion about the upper bound of  $\lambda(f)$ follows immediately from (iv).
\qed\endproof

Motivated by the BMY--inequality $c_1^2(X)\le 3c_2(X)$ (see \cite[Theorem VII, 4.1]{4authors}), one introduces  the \emph{proportionality deviation}    
\begin{equation}\label{eqn:BY}
\bmy(X):= 3c_2(X)- c_1^2(X).
\end{equation}
The BMY--inequality then becomes  $\beta(X)\ge0$ with  equality if and only if $X$ is a ball quotient.
\begin{corr}  \label{bmyformula} In the situation of \eqref{eqn:fibration}, one has
\[
\bmy(X)=4(b-1)(g-1) -\omega^2_{X/B} + 3\sum_{s\in B}  \delta_s.
\]
\end{corr}
The invariants $\delta_s$ are just the number of  double points of the fiber $X_s$. To show this, one uses
expressions for  the Euler number and the arithmetic genus of a semi-stable curve:
\begin{lemma} \label{BreakUpInvs} Let $C$ be a semi-stable curve on a compact complex surface $X$  and let the irreducible components
be  $C_\alpha$, $\alpha=1,\dots, N$; the genus of the normalization of $C_\alpha$  is denoted $g_\alpha$  and $\delta_C$ is  the number of double points of $C$. Then one has:\\
$$
\aligned
e(C)&= \sum_\alpha  (2-2g_\alpha)  -\delta_C.
\\
2p_a(C)-2&= K_X \cdot C+ C^2= \sum_\alpha  (2g_\alpha -2 ) +2\delta_C.
\endaligned
$$
\end{lemma}
\proof 1) follows from the additive property of the Euler number.\\
2) can be seen as follows. For an irreducible component $C_\alpha$ of $X_s$ one has $p_a(C_\alpha)=g_\alpha+\delta_\alpha$ with $\delta_\alpha$ the number of double points of $C_\alpha$; the number $\sum_{\alpha<\beta}  C_\alpha C_\beta$ is the number $\delta'_{s}$ of  those double points of $X_s$ that are intersections of two components. Hence
\begin{eqnarray*}
K_X\cdot C+ C^2& = &  \sum _\alpha (K_X\cdot C_\alpha + C_\alpha^2) + 2 \sum_{\alpha<\beta}  C_\alpha\cdot C_\beta \\
&= &\sum _\alpha (2p_a(C_\alpha)-2 ) + 2  \sum_{\alpha<\beta}  C_\alpha\cdot C_\beta \\
&=& \sum _\alpha (2g_\alpha  -2)+ 2(\underbrace{\sum \delta_\alpha + \delta'_s}_{\delta_C}). \qed 
\end{eqnarray*}
\endproof

\begin{corr}\label{stablefibresinvars} Suppose $f:X\to B$ is a semi-stable genus $g$ fibration. Then
\begin{eqnarray*}
\delta_s   =   \delta_{X_s}=\#( \text{\rm double points of } X_s).
\end{eqnarray*}
\end{corr}
\proof
 The previous lemma gives
 \begin{equation}\label{eqn:eufiber}
e(X_s)= \sum (2-2g_\alpha)  -\delta_{X_s}
\end{equation}
On the other hand, the arithmetic genus $p_a(X_s)$ is independent of $s\in B$ since by the adjunction formula  $2p_a(X_s)-2 = K_X\cdot X_s +X_s^2= K_X\cdot X_s$ which does not depend on $s$ since all fibers are numerically equivalent.  
So, applying the previous lemma once again, one finds
\begin{equation}
\label{eqn:eugen}
\begin{array}{lcl}
e(X_{\rm gen})=2-2g  &= & -(K_X\cdot X_s+ X_s^2)\\
 &=& \sum (2-2g_\alpha) - 2\delta_{X_s}.
\end{array}
\end{equation}
Comparing this with \eqref{eqn:eufiber} and using    \eqref{eqn:EgenEspec}, shows that 
$$\delta_s:= e(X_s)-e(X_{\rm gen})= \delta_{X_s} \qed .
$$
 \endproof

 \subsection{A refined BMY--inequality} \label{sec:BMY}
There is a refinement  of the BMY--inequality due to R. Kobayashi \cite{kob}, and Cheng-Yau \cite{cy}.  An algebraic proof has been given by  Miyaoka \cite{mi}, but that  proof gives no information in case of equality, information which is crucial
for the proof I give of Theorem (A).
\par
To state it, some preparations are needed. First, recall that a divisor $D$ on a surface is \emph{big} if  for $m\gg 0$ the linear system  $|mD|$  maps $X$ birationally onto its image. This is equivalent to $\kappa(D)=2$, where $\kappa$ denotes the Kodaira-dimension. 
If some multiple of $D$ is effective, $D$ is big precisely if  $D^2>0$ \cite[IV Prop. 7.4]{4authors}.
The divisor $D$ is  \emph{nef}, if  $D\cdot C\ge 0$ for all curves $C$. 

A surface $X$ by  definition is of general type if the canonical divisor $K_X$ is big, or, equivalently, if $\kappa(X)=\kappa(K_X)=2$.
A suitable multiple $|mK_X|$ thus  maps $X$ birationally onto its image. In this situation $K_X$ is nef precisely
when $X$ is minimal.  Even then, the bundle $K_X$  need not be ample:  $|mK_X|$  contracts all
  ADE--curve configurations. Such a configuration of curves  contracts to  a rational  double point under $|mK_X|$, $m\gg 0$. Any such rational double point   $p$ is a quotient singularity: its germ is of the form $U_p/G_p$  where $U_p$ is smooth and $G_p$ a finite group. 

The refinement gives an inequality, valid for  the so-called \emph{logarithmic Chern classes} for a  Zariski-open 
surface $X\setminus C$ where $C\subset X$ is a normal crossing curve. These are  defined as follows:
\[
\aligned 
c_2(X,C)&:=  c_2(X) -  e(C) ;\\
c_1^2(X,C)&:=   c_1^2(\Omega^1_X(\log C)) = c_1^2(X) +    2 K_X\cdot C + C^2.\\
\endaligned
\]
Introduce  the \emph{logarithmic proportionality deviation}
\[
\aligned
\bmy(X,C):=  & 3c_2(X,C)- c_1^2(X,C)\\
 =& \underbrace{3c_2(X) -c_1^2(X)}_{\beta(X)} -\underbrace{\left( 3e(C)+2K_X\cdot C+C^2\right)}_{\beta(C)} . 
\endaligned
\]
For every  A-D-E configuration  $R_p$, contracting to a singular point $p$, one introduces also  the following punctual proportionality deviations:
\begin{equation}\label{eqn:singterm}
\beta(p)= 3\left( e(R_p) - \frac{1}{|G_p|} \right) >  0.
\end{equation}
Now I can state the  promised strengthening of the BMY--inequality which results from Yau's techniques:
\begin{thm}[\protect{\cite[Thm 2]{kob}, \cite{mi,cy}}]  \label{kobaresult} Let $X$ be a compact complex surface,  $C\subset X$ a normal crossing divisor such that $K_X+C$ is nef and big.  Let $R_{p_j}$, $j=1,\dots, k$ be the A-D-E configurations in $X\setminus C$ and let $X'$ be the normal surface
 obtained after contracting  the $R_{p_j}$  to a singularity $p_j$.  Then
\begin{equation}\label{eqn:miyau}
\beta(X,C)- \sum_j \beta(p_j)  \ge 0
\end{equation}
and equality holds if and only if $X'\setminus C =  \Gamma  \backslash B^2 $, the quotient of the 2-ball $B^2$ by a discrete subgroup $\Gamma\subset \text{\rm PSU}(2,1)$ acting freely except over the singularities of $X'$ where $\Gamma$ acts with isolated fixed points. 
\end{thm}
\begin{rmk} \label{sscurves} If $K_X$ itself is nef, the result applies to the following two sorts of curves $C$ that are relevant here:
\begin{enumerate}
\item  any semi-stable curve on $X$ of arithmetic genus $\ge 2$;
\item  an elliptic curve  with   negative self intersection.
\item  disjoint unions of curves of  the above sort.
\end{enumerate}
To see that the condition of the  theorem  holds for such curves, argue as follows.
Let $C=F+R$ with $F$ a disjoint union of curves  as in 1, and $R$ a disjoint
union of curves as in 2.  If $F\not=\emptyset$, then by nefness of $K_X$, $(K_X+C)^2\ge (K_X+C)C =2 p_a(C)-2  >0$.
Otherwise,  $(K_X+R)^2=  K_X^2- R^2 >0$. 

Next, again using nefness of $K_X$, one has $(K_X+C)D \ge 0$ for all irreducible curves $D$  that are not
components of $C$ and thirdly, for all components $C_i$ of $C$ with $p_a(C_i)>0$, one has $(K_X+C)\cdot C_i \ge   K_XC_i+ C_i^2= 2p_a(C_i)-2\ge 0$. Finally,  a smooth rational  component $C_i$ of $C$ must meet the union $D_i$ of the other components in at least 2 points (this is a consequence of the  definition of semi-stablity) so that 
$$(K_X+C)\cdot C_i= K_X\cdot C_i +C_i^2+ D_i\cdot C_i \ge -2+ 2\ge 0.
$$
So $K_X+C$ is big indeed. See also \cite[Theorem 7.6.]{sak}.
\end{rmk}

 \subsection{Properties of excess invariants for curves}

 \begin{lemma}  \label{betaCurves}
(1) The excess invariant $\beta(C)$ is additive: 
for disjoint curves $C$ and $D$ one has  $\beta(C\coprod D)=\beta(C)+\beta(D)$.\\
(2) Let $C=\sum C_\alpha$ be a semi-stable curve with irreducible components $C_\alpha$ and total number $\delta_C$ of double points. Then, recalling that $g_\alpha$ is the genus of the normalization of $C_\alpha$, one has:
\[
\beta(C)=\sum_\alpha(2-2g_\alpha) +\delta_C -C^2.
\]
\end{lemma}
\proof 1) is an easy verification and 2) follows immediately from Lemma~\ref{BreakUpInvs}
since 
\[
\aligned
\beta(C)&= 3e(C)+2K_XC+C^2\\
 &= 3e(C)+ 2(K_XC+C^2)- C^2\\
 &= 3e(C)+ 4\sum_\alpha (g_\alpha-1)+4\delta_C- C^2\\
 &= 6 \sum_\alpha (1- g_\alpha) -3\delta_C + 4\sum_\alpha (g_\alpha-1)+4\delta_C -C^2. \qed 
\endaligned 
\]
\endproof

\begin{exmples}
1) Let $C=X_s$ a possibly singular fibre of a semi-stable fibration in curves of genus $g$.  Since $X_s^2=0$,
using  \eqref{eqn:eugen} in the formula for $\delta_s=\delta_{X_s}$ gives:
\begin{equation}\label{eqn:epsilon}
\beta(X_s)= 3\delta_s +2(1-g).
\end{equation}
This can indeed be positive and  in that case gives an amelioration of the BMY--inequality.\\
(2)  Consider   a fibre $X_s$ as in (2)  of the form $X_s=E_s+F_s+R$ with
$p_a(E_s)=1$, $F_s$  a     tree  of smooth rational curves meeting   $E_s$ and $R$ transversally in one point  and $R$,
 the union of the remaining components.  
 Note that if $X_s$ is   \emph{of compact type}, $E_s$ is necessarily
 a smooth elliptic curve and one can take for $F_s$ the tree of all rational curves (necessarily smooth) connecting $E_s$ with other curves of genus $>0$.
 While the curve $E_s$ is called an \emph{elliptic tail}, the divisor $F_s$ is  the corresponding  \emph{connecting rational tail}  of $X_s$.  By Remark~\ref{sscurves}  the semi-stable divisor $E_s$ can play the role of $C$ in Theorem~\ref{kobaresult}, but this is not the case
 for $E_s+F_s$ since  $K_X+E_s+F_s$ not nef: the last component of $F_s$ has negative intersection with $K_X+E_s+F_s$.
 Fortunately, the divisor  $F_s$ forms a chain of $(-2)$--curves and hence can be contracted to
 a rational double point, say $f_s$,  of type $A_{\ell_s}$.  
Make now note of the following inequalities
 \begin{equation}\label{eqn:epsilon2}
 \begin{cases}
 \beta(E_s) = -E_s^2 \ge 1=3\ell_s+1 & \text{if } F_s=\emptyset, \\
 \beta(f_s)= 3(\ell_s+1 - \frac{1}{\ell_s+1})\ge 3\ell_s+1 & \text{otherwise.}
  \end{cases} 
\end{equation}
where     Lemma~\ref{betaCurves} and \eqref{eqn:singterm}  have  been used.
   \end{exmples}

The preceding remarks have  the following useful consequence. Recall that a fiber $X_s$ is a singular curve of compact type
if and only if its dual graph is a tree and all components of $X_s$ are smooth.
\begin{corr}[of Theorem~\ref{kobaresult}] \label{omegaineq} Let  $X$ be a compact complex surface,
\[
f:X \to B 
\]
a semi-stable fibration of genus $g\ge 2$ curves. 
Let   $\Sigma\subset B$  be the   set  of critical  points of  $f$ over which the fiber is not of compact type.
Let $R$ be the disjoint union of the elliptic tails in fibers of compact type that
have no connecting rational tails and let $f_j$, $j\in J$ be  the rational double  points
obtained by contracting the  connecting rational tails. 
Setting  $C=R+\coprod _{s\in\Sigma} X_s$, assume that $K_X+C$ is nef and big.
%
Then one has 
\begin{equation}   \label{eqn:Ineq}
2(g-1)\deg \Omega^1_B(\log (\Sigma)) -\omega_{X/B}^2 + 
  3\sum_{s\in \Sigma_c} \delta_s -\beta(R) -\sum_{j\in J}  \bmy(f_j) \, \ge \, 0.
\end{equation}
\end{corr}
\proof Let   me calculate $\bmy(X,C)$. 
Additivity of $\beta$ and \eqref{eqn:epsilon} give
\[
\bmy(C)= 3 \sum_{s\in \Sigma} \delta_s +\beta(R)+2(1-g)\cdot \#\Sigma .
\]
On the other hand, by   Cor.~\ref{bmyformula}, setting  $ b=\text{\rm genus}(B)$, one has
\[
\aligned
\bmy(X)&=  4(b-1)(g-1) -\omega_{X/B}^2 +3\sum_{s\in B} \delta_s\\
 &= 2(g-1)\deg \Omega^1_B(\log (\Sigma)) -\omega_{X/B}^2 +2(1-g)\cdot\#\Sigma + 3\sum_{s\in B} \delta_s   \\
\endaligned
\]
and so
\[
\aligned
\bmy(X,C) &= \bmy(X) -\bmy(C)\\
 &=  2(g-1)\deg \Omega^1_B(\log (\Sigma)) -\omega_{X/B}^2 + 3\sum_{s\in B} \delta_s - 3\sum_{s\in \Sigma} \delta_s
 -\beta(R)\\
 &= 2(g-1)\deg \Omega^1_B(\log (\Sigma)) -\omega_{X/B}^2 + 3\sum_{s\in \Sigma_c} \delta_s -\beta(R).
\endaligned
\]
So the desired   inequality follows from   Theorem~\ref{kobaresult}. \qed
\endproof
\begin{rmk} The inequality \eqref{eqn:Ineq} will be used in the following form:
\begin{equation}
\label{eqn:Ineq2}
\left. \aligned
 4\left(1-\frac 1g\right) +  \frac{1}{d(f)} \cdot  \left( \sum_{s\in \Sigma_c}  3\delta_s  -\beta(R)-\sum_{j\in J} \beta(f_j)\right)
 &  \ge \frac{\omega_{X/B}^2}{d(f)} \\
  d(f) := \frac{g}{2}  \deg \Omega^1_B(\log (\Sigma)).&
\endaligned\right\}
\end{equation}
\end{rmk}

 \subsection{On Moriwaki's slope inequality}\label{sec:Moriwaki}
 Xiao's slope inequality \eqref{eqn:xiao} admits a refined version  which is due to  Moriwaki. It is a consequence of  certain enumerative properties of cycles on the compactified moduli space $\overline \MM_g$. To explain it, let me introduce some more notation. The moduli space $\overline \MM_g$ has as boundary the irreducible divisors $\Delta_j$, $j=0,\dots,[\frac g 2]$ where the stable genus $g$ curve $C$ belongs to  $\Delta_0$ if $C$ is irreducible and to $\Delta_j$, $j>0$, if it is of the form $C=C_1+C_2$ where $p_a(C_1)=j$, $p_a(C_2)=g-j$. The double point $C_1\cap C_2$ is then  called a double point of type $j$. The double points of an irreducible curve are called  of type $0$.  This terminology  extends  to the semi-stable situation: 
 if there is a tree $R$ of rational  curves connecting   $C_1$ and $C_2$, all of the double points on $R$ are of type $j$.

If $p : B \to \overline{\MM_g}$ is the period map for $f$, the curve $p(B)$ meets the divisor $\Delta_j$ exactly in the points $s\in B$ over which there is a singular fibre  with a double point of type $j$. Let $\delta_j(f)$ be the total number of such points. Moriwaki's inequality  \cite[Theorem D]{moriwaki} reads:
\begin{equation}\label{eqn:lastineq}
(8g+4) \deg f_*\omega_{X/B} \ge g \cdot \delta_0(f)+ 4 \sum_{j=1}^ {[\frac g 2]}  j(g-j) \cdot  \delta_j(f).
\end{equation}
To see that this is a refinement of the slope inequality
\eqref{eqn:xiao}, observe that by Lemma~\ref{stablefibresinvars}  $\sum_{s\in B} \delta_s$ is the total number of double points $=\sum_{j=0}^{[\frac g 2]} \delta_j(f)$ and so the right hand side of \eqref{eqn:lastineq}
reads
\begin{eqnarray*}
\begin{array}{clc}
   g (\sum_{s\in B} \delta_s) &+ & \hspace{-12em}  \sum_{j=1}^ {[\frac g 2]} \left( 4 j(g-j)-g  \right) \cdot \delta_j(f)  \\
  &=& g (12 \deg f_*\omega_{X/B} -\omega^2_{X/B}) +    \sum_{j=1}^ {[\frac g 2]} \left( 4 j(g-j)-g \right) \cdot\delta_j(f).
\end{array}
\end{eqnarray*}
Note that the second line uses Lemma~\ref{basics}.(iv).
Dividing   the   inequality  \eqref{eqn:lastineq}  by $g\cdot \deg f_*\omega_{X/B} $  then indeed  leads to the sharpening  of \eqref{eqn:xiao}:
 \begin{equation}\label{eqn:moriwaki}
\lambda(f) \ge 4- \frac 4g + \frac{1}{\deg f_*\omega_{X/B}}\left[ \sum_{j=1}^ {[\frac g 2]} \left( 4 \frac{ j(g-j)}{g}-1 \right)\cdot\delta_j(f) \right]
\end{equation}

 \section{Background from Hodge theory}
 \label{sec:Hodge}
 
 \subsection{Higgs bundles and the Arakelov inequality}
  \label{sec:HiggsAndArak}
For background on Hodge theory used in this section see for example
\cite[Chapter 13.1,13.2]{periodbook} and \cite[Chapter 11]{mhs}.

The central concept in this section is that of a Higgs bundle:
\begin{dfn} A (logarithmic) \emph{Higgs bundle}  on a curve $B$ with  poles in  a finite set $\Sigma\subset B$ is a vector bundle $\HH$ on $B$ 
with a sheaf morphism $\sigma : \HH\to \HH\otimes \Omega^1_B(\log \Sigma)$, the \emph{Higgs field}.\footnote{Higgs fields $\sigma$ over a higher dimensional base also have to satisfy $\sigma\wedge \sigma=0$} A \emph{graded   Higgs bundle}
$\HH=\oplus \HH^j$ is a Higgs bundle verifying verify $\sigma|\HH^j \to \HH^{j-1}\otimes  \Omega^1_B(\log \Sigma)$.
\end{dfn}
A variation  of Hodge structure naturally gives a  graded Higgs bundle with Higgs field induced by the Gauss-Manin connection.
Let me explain this in the current situation. So assume that one has a polarized variation of Hodge structure of weight $1$ over 
$B_0=B\setminus \Sigma$ with underlying  local system $\mathbf{H}$. 
The local monodromy of $\mathbf{H}$ around points of $\Sigma$ is assumed to be unipotent.\footnote{In the geometric situation
this is always the case since the degenerating fibers are semi-stable.}
The vector bundle $\mathbf{H}\otimes_\C \OO_{B^0}$ together with the Gauss-Manin connection then 
turns out to have a canonical extension as a vector bundle $\HH$ on all of $B$ equipped with a  connection  with 
logarithmic poles in points of $\Sigma$, say
\[
\nabla: \HH \to \HH \otimes_\C \Omega^1_S(\log \Sigma). 
\]
The Hodge bundle $\FF^1\subset \mathbf{H} \otimes_\C \OO_{B^0}$ admits a unique   extension 
$$
\HH^{1,0} \subset \HH
$$
and hence a two-step filtration $\HH^{1,0}\subset \HH$. On the associated graded  bundle $\HH_{\rm Hgs}=\HH^{1,0}\oplus \HH^{0,1}$
 the logarithmic  connection  $\nabla$ defines the structure of a Higgs bundle with
 (logarithmic) Higgs field 
\[
\sigma  :  \HH^{1,0}  \to \HH^{0,1}\otimes \Omega^1_B(\log \Sigma).
\]
A  graded Higgs subbundle of $\HH_{\rm Hgs}$ consists of a graded holomorphic subbundle preserved by $\sigma$. If such a subbundle
comes from  a local subsystem of $\C$-vector spaces, and the Higgs field vanishes on it, the  polarization 
makes the  local system into  a \emph{
unitary local system}.\footnote{If this system would be defined over $\Q$ its monodromy would be finite and the Higgs subbundle corresponds to
an isotrivial subsystem.} On the other side of the spectrum we may have Higgs subbundles for which the Higgs field is an isomorphism.  

A central property is the semi-stability of Higgs bundles:
\begin{prop}  \label{SimpsonCorr}The first Chern form of a graded Higgs subbundle 
$\GG$ of  a   Higgs bundle $\HH$  coming from a \vhs\ on $B\setminus \Sigma$ (with respect to the Hodge metric) is negative
semidefinite, and it is $0$ everywhere if and only if $\GG^\perp$ is a
(holomorphic) graded Higgs subbundle as well.  In this case the
 variation splits as  an orthogonal direct sum of Higgs bundles $\HH=\GG\oplus \GG^\perp$.  
\end{prop}
For an elementary proof see \cite[Chapter 13.1]{periodbook}. 
For the weight one case, the  above corollary  applies in particular to  the maximal unitary subsystem   $\HH_{\rm un}$
on which $\sigma=0$. Its Chern form is identically zero and hence  by Prop.~\ref{SimpsonCorr} there is  an orthogonal sum  decomposition
\begin{equation}
\label{eqn:Splitting}
(\HH_{\rm Hgs},\sigma)= (\HH_{\rm max},\sigma)\oplus (\HH_{\rm un},0).
 \end{equation}
where the first Higgs bundle has injective Higgs field.

Semi-stability  implies the following Arakelov type inequality:

\begin{thm}[ \cite{faltings}] \label{arakelov}
Let $(\HH_{\rm Hgs},\sigma)$ be  a  graded Higgs bundle over the curve $B$ with poles in $\Sigma$ which
induced by a polarized weight one variation of Hodge structure as sketched above.
The  following inequalities for  the Higgs component  $\HH^{1,0}$ hold:
\begin{equation} \label{eqn:Arak0}
0   \le  \deg \HH^{1,0}  \le 
 \frac 1 2 \cdot     \rank \HH^{1,0} \cdot  \deg \Omega^1_{X/B}(\log \Sigma) .
\end{equation}
Moreover, if      the Arakelov bound for $\HH^{1,0}$  is attained,   then
$\HH_{\rm Hgs}$ has   no unitary  subsystem, i.e. the splitting \eqref{eqn:Splitting} just reads
$\HH_{\rm Hgs}=\HH_{\rm max}$,  and, moreover,   the Higgs field is an isomorphism. 

Conversely, if  the Higgs field for
$\HH_{\rm Hgs}$ is an isomorphism, the Arakelov bound for $\HH^{1,0}$ is attained. 
 \end{thm}
 \proof
  The first inequality is easy: $(\HH^{1,0},0)$ is a quotient Higgs bundle of $(\HH_{\rm Hgs},\sigma)$ and hence, by semi-stability (Theorem~\ref{SimpsonCorr}) has degree $\ge 0$.

 For the second inequality, note that   the degree of the components of the Higgs bundle $\HH_{\rm Hgs}$ and $\HH_{\rm max}$ are
 the same while  the Arakelov inequality also holds for $\HH_{\rm max}$. So, if the Arakelov inequality is attained for $\HH_{\rm Hgs}$, one must have $\HH_{\rm Hgs}=\HH_{\rm max}$. I shall assume this from here on. 
Consider the  Higgs  bundle saturation  $\GG=\HH^{1,0}\oplus \GG^{0,1}$ for $\HH^{1,0}$, i.e. we consider
 \[
 \sigma: \HH^{1,0}\onto \GG^{0,1}\otimes \Omega^1_B(\log \Sigma),
 \]
 where by definition of the saturation  the right hand side is the image of the Higgs field. The right side has degree
 $
(  \rank \GG^{0,1}) \cdot\deg \Omega^1_{X/B}(\log \Sigma)
 $
and  since the pair $(\ker\sigma ,0)$ is obviously a Higgs subbundle of $(\HH_{\rm Hgs},\sigma)$, by Cor.~\ref{SimpsonCorr}, its degree is non-positive which implies
 \begin{equation}
 \label{eqn:TrivIneq1}
 \deg \HH^{1,0}  \le \deg \GG^{0,1}+ (\rank \GG^{0,1} ) \cdot\deg \Omega^1_{X/B}(\log \Sigma).
 \end{equation}  
On the other hand, by stability   $\deg\GG= \deg \HH^{1,0}+ \deg \GG^{0,1} \le 0$ and hence, adding these two ineqalities, 
the second inequality follows. Moreover, if equality holds,   $\deg(\GG)=0$  and $\GG=\HH_{\rm Hgs}$ so that $\sigma$ is surjective.
Since $\HH_{\rm un}=0$, the Higgs field has no kernel and hence $\sigma$ is an isomorphism.
The proof of the converse statement is left to the reader.\qed
 \endproof
\begin{corr} \label{nontrivialfixedpart}
Let $f: X\to B$ be a  family of semi-stable curves of genus $g$ over a curve $B$ with $\Sigma\subset B$  the subset of 
those critical points $s\in B$ of $f$  for which  $\jac(X_s)$ is \emph{not}  compact.  For the associated logarithmic Higgs bundle 
$(\HH_f,\sigma_f)$  let $\HH_f=\HH_{f,\rm max} \oplus \HH_{f,\rm un}$ be the splitting \eqref{eqn:Splitting}. With    $g_0$
the rank of the bundle     $\HH_{f,\rm un}^{1,0}$ ,  the Arakelov 
inequality refines to
\begin{equation}
  \deg \HH^{1,0}_f =\deg \HH^{1,0}_{f,\rm max}   \le 
 \frac 1 2 \cdot     (g-g_0) \cdot \deg \Omega^1_{X/B}(\log \Sigma) . \label{eqn:Arak}
\end{equation}
The Arakelov bound  for the full Hodge bundle $\HH^{1,0}_f$ is attained, i.e.,
\begin{equation}
\label{eqn:ArBndAtt}
\deg f_*\omega_{X/B} =\half g  \cdot\deg  \Omega^1_B(\log \Sigma)
\end{equation}
 if and only if $\HH_f=\HH_{f,{\rm max}}$ and  the Higgs field
$ \sigma_f$  is an isomorphism. 
\end{corr}
\begin{rmk}  The proof shows that  the Arakelov bound for $\HH_{f,\rm max}^{1,0}$ is attained if and only if its own Higgs field is an isomorphism.
\end{rmk}

Below I also need a rigidity result for Higgs bundles.  Recall that a  deformation of a family $f:Y  \to B$,  smooth over $B_0=B\setminus \Sigma$  and  fixing  $(B,\Sigma)$, is a morphism
$F:(Y,E)\to  (B,\Sigma)\times T$ with $T=(T,o)$ a germ of a complex manifold, such that $F|(B,\Sigma)\times \set{o}=f $. 
The  family $f$ is called \emph{rigid} if  every deformation $F$ of $f$ is induced from $f$ by pullback along the projection $(B,\Sigma)\times T\to (B,\Sigma)$. 
 Since for a family $Y$ of Abelian varieties
 the  infinitesimal Torelli theorem holds, assuming $f$ is not locally trivial over $B_0$, rigidity for $f$
is equivalent to rigidity for the associated period map $p(f): B_0 \to \AA_g$.
\begin{prop}[\protect{\cite[Prop. 3.7.]{arakelovstuff}}] \label{Rigidity} Assume that $f$ is not locally trivial over $B_0$ 
or, equivalently, that the period map   $p(f)$ is not constant. If  the associated Higgs field  for $\HH_f$ is an isomorphism, 
then $p(f) $ is rigid. 
In particular,   if $f$ reaches the Arakelov bound for $\HH_f^{1,0}$,  the curve $B_0$ is 
rigidly embedded in $\AA_g$.
\end{prop}

 \subsection{Application of the slope inequality to   fibrations attaining the Arakelov bound}
 \label{sec:mainresult}

\begin{thm} \label{mainthm} For a  non-isotrivial   genus $g\ge 2 $ fibration $f:X\to B$ of semi-stable  curves 
attaining the Arakelov bound,  one has $g=2,3,4$.
\end{thm}
\proof  Step 1. One first  has to make sure that the condition of Theorem~\ref{kobaresult} and Corollary~\ref{omegaineq} are verified, i.e., that  $K_X+C$ is nef and big for
$C=R+\coprod _{s\in\Sigma} X_s$, $R$ the union of certain elliptic tails.
If   $B$ is not rational, the surface $X$ is minimal: no multisection of $f$ can be rational and $f$ may be  assumed to be relatively
minimal. Remark \ref{sscurves} shows that it suffices to prove that $K_X$ is nef, or equivalently  \cite[VI. Theorem 2.1]{4authors} that
the Kodaira dimension of $X$ is $\ge 0$. Since $g\ge 2$, if $b\ge 1$  the Iitaka-inequality \cite[III,Theorem 18.4]{4authors} shows
that  in fact $\kappa(X)\ge 1$.

The case $b=0$ remains. Since $f$ is not isotrivial , by \cite[III,Theorem 18.2]{4authors}, one has $\deg ( f_*\omega_{X/B})>0$.
Since one  may also assume that $g\ge 5$,  the Arakelov equality  then gives  $\deg ( f_*\omega_{X/B})= \half g (-2 +\# \Sigma)\ge  3$.
Riemann-Roch for $f_*\omega_X^{\otimes n}$  thus  yields the following estimate for the  plurigenera
\[
P_n(X) \ge H^0(B,f_*\omega_{X/B}^{\otimes n})\ge   n\deg(f_*\omega_{X/B}) -2n + 1\ge n+1 
\]
and  hence  $\kappa(X)\ge 1$.    If $X$ is minimal, this finishes the proof.
However,  $X$ could fail to be minimal.   This happens if  there exists  exceptional curves $E$ 
which figures as a (multi)-section.  

Write $K_X= K_{\rm min}+ \sum E$, where $K_{\rm min}=\sigma^* K_Y$, $\sigma: X\to Y$ the contraction of all the minimal curves $E$ to a minimal model $Y$.
Then $K_{\rm min}$ is nef on $X$ and so, by the  argument of Remark~\ref{sscurves},
$(K_X+C)\cdot D\ge 0$ for all curves $D$ that are not exceptional.  
Furthermore, since $|\Sigma|\ge 3$, $E$ being a (multi)section, one has $C \cdot E\ge 3$ and since $K_X\cdot E=-1$,
it follows that  $(K_X+C)\cdot E  \ge 2$. Hence $K_X+C$ is nef. To show that it is also big, one estimates its self-intersection
as follows. Using nefness of $K_X+C$ and the previous estimate $(K_X+C)\cdot E  \ge 2$, one finds
\[
(K_X+C)^2\ge (K_X+C)K_X= (K_X+C)\cdot (K_{\rm min}+\sum E)  \ge  2\# \text{(exc. curves } E) >0,
\]
since $X$ is not minimal.
\medskip
 
\noindent Step 2:  \emph{Kodaira fibrations cannot attain the  Arakelov bound}. Indeed,  $\omega_{X/B}^2=12\deg f_*\omega_{X/B}=12g(b-1)$ while, by Cor.~\ref{omegaineq}  $\omega_{X/B}^2\le 4(g-1)(b-1)< 12g(b-1)$. It follows that at least one $\delta_i(f)$ is non-zero.
 
 \medskip
\noindent  Step 3: \emph{The case where $\delta_j(f)>0$ for some $j>0$}.  
As remarked in Step 1,   $\deg f_*\omega_{X/B}>0$ since  $f$ is not locally trivial and so
the Moriwaki slope inequality \eqref{eqn:moriwaki}  applies. On the other hand,
 in  the inequality  \eqref{eqn:Ineq2}, by the assumption on the Arakelov bound, one has $d(f)= \deg  ( f_*\omega_{X/B})$,
and so this inequality combined with   Moriwaki's slope inequality gives:
\begin{equation*} 
\aligned
 4 - \frac  4 g   +   \frac{1}{ \deg f_*\omega_{X/B}} 
    & \left [\sum_{s\in  \Sigma_c}  3\delta_s -\beta(R)-\sum \beta(f_j)\right]   \\
    \ge  \,   \lambda(f) &=   \frac{ \omega^2_{X/B}}{d(f)} \\
\ge \,  &4- \frac 4g +     \frac{1}{\deg f_*\omega_{X/B}}   \left[ \sum_{j=1}^ {[\frac g 2]} \left( 4 \frac{ j(g-j)}{g}-1 \right)\cdot\delta_j(f) \right].
%
%
\endaligned 
\end{equation*}
One gets
\begin{equation}\label{eqn:curxineq}
\sum_{j=1}^ {[\frac g 2]} \left( \frac{ 4  j(g-j)}{g}-1 \right)\delta_j(f)  -3\sum_{s\in \Sigma_c} \delta_s  + \beta(R)+\sum_{j\in J} \beta(f_j)  \le 0.
\end{equation}
Moreover, if the left  hand side is zero, then $\beta(X,C)+\sum_{j\in J} \beta(f_j)=0$ where
$
C=R+\coprod _{s\in\Sigma} X_s
$.

Rewrite \eqref{eqn:curxineq}  as a sum of terms each  involving double points of fixed type.
First consider the contribution of double points of type $1$.  If such a double point   belongs to a fiber of non-compact type,
 it contributes $4(1 -\frac 1 g )-1=3 -\frac 4g$.  On  a fiber  $X_s$ of compact type  these are the double points on the 
  semi-stable  elliptic tail
 $E_s+F_s$ plus one further point where the tail connects to the remaining components of the fiber, and so there are $\ell_s+1$  double points of type $1$ in $X_s$.  Now recall inequality \eqref{eqn:epsilon2}: 
whether or not a connecting rational tail exists or not, the contribution of such an elliptic tail  $E_s$ to the third and fourth term in  \eqref{eqn:curxineq}   is at least   $3\ell_s+1$.
 It follows that  the total contribution
 to \eqref{eqn:curxineq} 
  is  at least
 $$
\underbrace{ (\ell_s+1)( 3 -\frac 4g )}_{1\text{st term}} \underbrace{-3( \ell_s+1)}_{2\text{nd term}}  + 3\ell_s+1\ge( \ell_s+1) \left(1 - \frac 4g\right).
 $$
 So for both situations  the "average" contribution  from  each  double point of   type  $1$ is 
  at least $1 - \frac 4 g$.  
  
  For  double points of type $j\ge 2$, the first term contributes at least $7 -\frac {16}{g}$ and   the second term contributes  $0$ or   $-3$
  and so their contribution is at least $4 -\frac {16}{g}$.
Hence, in total, one gets:
\begin{eqnarray*}
0&\ge & \sum_{j=1}^ {[\frac g 2]} \left( \frac{ 4  j(g-j)}{g}-1 \right)\delta_j(f)  -3\sum_{s\in \Sigma_c} \delta_s  +\beta(R)\\
&\ge &(1 - \frac 4 g )\cdot \delta_1(f) +  \left[ 4 -\frac {16}{g} \right]\cdot  \sum_{j \ge 2} \delta_j(f)
\end{eqnarray*}
 It follows that $g\le 4$. 
 \medskip

 \noindent Step 4: \emph{The case  $\delta_j(f)=0$ for  all $j>0$}. I shall argue, following a  suggestion
 of Kang Zuo,  that this case does not occur, thereby completing the proof. So suppose it  does occur.
Note that  in that case  there are only double points of type $0$. These might include the double points
in a string of $(-2)$ curves. However,  \eqref{eqn:curxineq} is an equality and since there are no
elliptic tails since $\delta_1(f)=0$, the above  argument implies that  $\beta(X,C)=0$. So, by Theorem~\ref{kobaresult}  the open surface  $X^0=X\setminus C$ does not  contain $(-2)$ curves. Furthermore, $X^0$  is a locally symmetric space for the group $\text{\rm PSU}(2,1)$ which is of  rank $1$. So  its Baily-Borel compactification
$X^*$  consists of $X\setminus C$  together with some zero-dimensional boundary components. Since $X^*$ is the minimal compactification, there is a holomorphic map $X\to X^*$ which necessarily must contract the curves  from $C$, by assumption 
a non-empty union of semi-stable fibers $X_s$, each consisting of an irreducible curve  of arithmetic genus $g$ plus, possibly some rational $(-2)$-curves.  As for any fiber configuration, the intersection matrix of the components of $X_s$  is negative
semi-definite with one dimensional kernel. This is  impossible by Mumford's   contraction criterion \cite{blowdown} which states
that the intersection matrix of a blown down configuration of curves is strictly negative definite.
This contradiction shows that this case does not occur.\qed
 \endproof

 \section{On the Coleman-Oort conjecture}\label{sec:COconj}
\subsection{Special subvarieties of Shimura varieties} \label{sec:SShimura}

For this subsection the reader is advised to consult \cite{Torloc}.
\medskip
\\
To start with, it is helpful to consider bounded symmetric domains $D=G(\R)/K$ from the  perspective of Shimura domains
as follows. 
Assume for simplicity that  $G(\R)$ is the group of holomorphic automorphisms of $D$. This is  the group of real points
of  a connected semi-simple  $\Q$--algebraic  group $G$
of adjoint type. It acts transitively on $D$ and  $K=K_x$ is
the isotropy group of some point $x\in D$. One shows that it comes with a  homomorphism 
\begin{equation}\label{eqn:HdgHom}
h_x: \bS \to G(\R) , \quad  \bS= \Res_{\C/\R} \C^\times,
\end{equation}
where "Res" stands for the \emph{Weil restriction}; here this simply means  that $\bS$ is   the group $\C^\times$ considered as a real algebraic group. 
The group $G$ acts on   $h_x$ by 
conjugation and one has  a natural identification
\[
D =\set{ \text{conjugacy class of } h_x}; \quad g\cdot x \iff h_{g\cdot x}= g h_x g^{-1}.
\]
This is the case, because the centralizer of   $h_x$  is the center of $K_x$ so that $h_{g\cdot x} = h_{\tilde g \cdot x}$ if and only if $g\tilde g^{-1}$ belongs to this center and hence the conjugate by $g$ of $h_x$ only depends on the 
class of $g$ modulo $K_x$. 

This description of $D$ is now indeed  essentially the description of a \emph{Shimura domain} or a \emph{Shimura datum}: a pair $(G,D)$
with $G$ a reductive $\Q$--algebraic group and $D$ a conjugacy class of  morphisms  \eqref{eqn:HdgHom} satisfying  some extra technical
conditions \cite[(2.1.1.--3)]{shimuravars}, which are not important here.
For any  arithmetic subgroup $\Gamma$ of $G$, the quotient $\Gamma\backslash D$ is called a \emph{Shimura variety}.

Let  $(D',G')$ be another Shimura domain and let there be given  a  homomorphism  $f:G \to G'$ of $\Q$--algebraic groups. Then $f\comp h_x$ defines a point $\bar f(x)\in D'$ and   $(f,\bar f): (G, D)  \to (G',D')$ is  an equivariant holomorphic map between bounded symmetric domains. 
For simplicity, let me assume that $f$ and hence $\bar f$ are injective; we thus
 we may identify $D$ with its image in $D'$. 
 One then says that $D$ is a \emph{Shimura subdomain} of $D'$ and if $f(\Gamma)\subset \Gamma'$ for some arithmetic subgroups
$\Gamma\subset G$, $\Gamma' \subset G'$ respectively, one says  that $\Gamma\backslash D$ 
is a \emph{special subvariety} of the Shimura variety $\Gamma'\backslash D'$. 

%

\subsection{Special curves in $\AA_g$}\label{sec:ShimCurves}

By Satake's classification \cite{satake}, there are  two essentially different classes of Shimura varieties of dimension $1$, also
called \emph{Shimura curves}:  

\begin{itemize}
\item Non-compact Shimura curves. These  are precisely the  \emph{modular curves}:   quotients of a Shimura domain  $(\SL 2 ,\eh)$ by a modular subgroup  $\Gamma\subset \SL{2;\Q}$. After replacing $\Gamma$
by a suitable finite index subgroup, these  carry a  universal family $E_\Gamma \to C_\Gamma=\Gamma\backslash \eh$ of elliptic curves. Here $\eh$ denotes the upper half plane.

\item   Compact Shimura curves which all arise as follows.   Let $F$ be a totally real number field and $\mathbf D$ a quaternion algebra
over $F$ which splits at exactly one place. Consider 
\[
\su {1,\mathbf D} =\sett{ x \in \mathbf D^\times}  {{\rm Nrd}_{\mathbf D}(x)=1} 
\]
 and the arithmetic subgroup
 \[
\Gamma_{\mathbf D} =\sett{ x\in \text{ some fixed maximal order of }\mathbf D}{ {\rm Nrd}_{\mathbf D}(x)=1}.
\]
Then, $(\Res_{F/\Q}\su {1,\mathbf D} , \eh)$ is a Shimura datum\footnote{As before "Res" stands for the Weil restriction; here
it means that the  $F$--group in question is to be considered as a group over the field $\Q$.}
and 
\[
C_{\Gamma_{\mathbf D}}:=  \Gamma_{\mathbf D}\backslash \eh
\]
is  a compact Shimura curve. Note that with $\bH$ the algebra of the quaternions, one has
\[
\aligned
\mathbf D(\R) &\simeq M(2;\R) \otimes_\R  \underbrace{\bH \otimes_\R\cdots \otimes_\R \bH}_{m-1},\quad m=[F:\Q] \\
G(\R) & \simeq \SL{2,\R}\times \underbrace{\su 2 \times\cdots \times \su 2 }_{m-1}.
\endaligned
\]
which shows that indeed, the upper half plane is the associated bounded domain.
Again, after possibly replacing $\Gamma_{\mathbf D}$ by a subgroup of finite index, there is a universal family of abelian varieties $A_{\Gamma_{\mathbf D}} \to C_{\Gamma_{\mathbf D}}$ over this Shimura curve.
The fibers of this family turn out to have
  dimension  $g=2^{m-1}$ or $g=2^m$.
  \end{itemize}

In order to classify special curves in $\AA_g$, one needs to understand  the embeddings of these curves  
 in $\AA_g$ as special subvarieties.
There are many  non-rigidly embedded Shimura curves; this depends  on the arithmetic of the corresponding embedding
of $\Q$--algebraic groups $G \to \smpl g$,  where $G$ is the   $\Q$-simple group of the type  introduced above (i.e., $G(\R)$ has precisely one non-compact factor, $\SL{2;\R}$).
In view of  Prop.~\ref{Rigidity}, I shall only consider rigidly embedded Shimura curves.  The result is as follows:

\begin{prop} \label{RigidSHiCurves}  Each of the two types of Shimura curves $C_\Gamma$, respectively  $C_{\Gamma_D}$ can be rigidly embedded in  $\AA_g$ for some $g$.  For $k$ large enough,  one also gets a  rigid embedding   
   of $C_\Gamma$, respectively $C_{\Gamma_D}$   in the  fine moduli space  $\AA_{g,k}$  of polarized Abelian varieties of dimension $g$ with level $k$--structure. 
  \\
  1. For the first type, $g$ can be arbitrary and  the modular family  of Abelian varieties is isogeneous to  a fibered  product 
$$
\underbrace{ E_\Gamma \times_{C_\Gamma}  \cdots \times _{C_\Gamma}  E_\Gamma }_{g \text { copies. }}.
$$
2.  For the second type, $g$ is a  multiple of 
$g=2^m$ or $2^{m-1}$, $m=[F:\Q]$ as before and
the modular family is isogenous to a fibered product
$$
A_{\Gamma_D} \times_{C_{\Gamma_D}}   \cdots \times_{C_{\Gamma_D}}  A_{\Gamma_D}
$$
of the modular family over a Mumford curve. 
\end{prop}
\noindent\textit{Sketch of proof.} If follows from  \cite[\S 6]{Saitoiii} that  
rigid representations in the sense of  \cite{addington}  give precisely the rigid embeddings we are classifying. 

In loc. cit. is explained that symplectic  representations are all constructed as follows.  As before, let  $F$ be a totally real number field and ${\mathbf D}$ a quaternion 
algebra over $F$ which splits at exactly one place.
For simplicity, assume that $F/\Q$ is Galois with Galois group $\GG$.
 The set of real embeddings $i_\alpha:F\into \R$, $\alpha=1,\dots,m$ 
is identified with $S=\set{1,\dots,m}$, the unique non-compact place being the first. We have 
\[
G^\alpha = \su{1,{\mathbf D}} \otimes_{i_\alpha} \R,\quad  G^\alpha(\C)=\SL{V^\alpha,\C},\quad V^\alpha=\C^2.
\]
The representations $V^\alpha$ are called atoms; molecules are the  tensor representation $V^{\alpha_1}\otimes\cdots \otimes V^{\alpha_s}$,  $\alpha_1,\dots,\alpha_s\in S$ 
and a direct sum of tensor representations is called a polymer. The Galois group $\GG$  acts on polymers. Polymers giving
rigid representations 
are  $\GG$--stable polymers consisting of a special type of molecules:  If $m=1$ only one atom per molecule is allowed and
if $m>1$ each molecule contains exactly one "non-compact" atom, i.e.,   $V^1$.  Note that  if $m=1$ the only such polymers are direct sums of the atom $V^1$ and in the 
second case one  can only have   direct sums of the 
molecule $V^1\otimes\cdots\otimes V^m$.  If $m=1$ these representations are defined over $\Q$ and below  
a symplectic embedding is exhibited. For $m>1$ some more work is need.
\\
Case 1. Non-compact embedded curves. Satake in \cite{satake} classifies  embeddings of algebraic  groups $G \into \smpl g$ for
  $G$  a  $\Q$--simple algebraic group
  and $G(\R)$   without compact factors. For $G=\SL 2$ such a  Satake embedding is induced by the following symplectic embedding.
Set $V_k=(\R^{2k},J_k)$, $J_k=\begin{pmatrix}
0 & \mathbf{1}_k\\
-\mathbf{1}_k & 0
\end{pmatrix}.$  The direct sum $\oplus_k V_1$ is isomorphic to the symplectic space  $V_k$. Whence  a  faithful representation {
$\rho_k$ of $\SL 2 $:
\[
\begin{pmatrix}
a& b\\
c& d
\end{pmatrix} \mapright{\rho_k} \begin{pmatrix}
a\mathbf{1}_k & b\mathbf{1}_k\\
c \mathbf{1}_k& d\mathbf{1}_k
\end{pmatrix}.
\]
For any $k=1,\dots,g$ the direct sum representation $\rho_k\oplus $ (rank $(g-k)$ trivial representation)
  induces a holomorphic embedding $\eh \into \eh_g$. It  gives the non-compact  embedded Shimura  curves starting from the Shimura datum $(\SL 2 ,\eh)$.  
There is no locally constant factor if and only if $k=g$ and then  the embedding is rigid. 
 The non-compact rigid curves are often called  \emph{rigid curves of Satake type}.  From  its  description one sees that the 
universal family is the fibered self product of the universal family of elliptic curves over the upper half plane.\\
Case 2. Symplectic embeddings when compact factors are present. 
Let ${\mathbf D}^\alpha$ be the twisted $F$--algebra $D$  obtained via the embedding $i_\alpha$
and form the  "corestriction from $F$ to $\Q$ of the algebra ${\mathbf D}$":
\[
\mathbf E= \text{Cor}_{F/\Q} {\mathbf D} :=  {\mathbf D}^1 \otimes_F \cdots  \otimes_F {\mathbf D}^m.
 \]
There are two possibilities for this algebra:
\[
\mathbf E=\begin{cases} M( 2^m,\Q) & \text{(split case})\\
M(2^{m-1}, H),\, H \text{  a quaternion algebra}/\Q & \text{(non-split case)}.
\end{cases}
\]
The $G$--representation space
\[
W=   V_1\otimes V_2\otimes\cdots\otimes V_m,\quad V_\alpha=\C^2\text{ with $G^\alpha(\C)$ acting }
\]
turns out to  be a  symplectic representation over $\Q$ in the split case ; in the non-split  case this is true for $W\oplus W$. See \cite[Sect. 4]{addington} or   \cite[Section 5]{carshcurves}.  The corresponding embedding
$$
\eh_1 \into \eh_g,\quad g=2^{m-1}, \text{ resp. } g=2^{m}
$$
is then a rigid embedding, as explained above. It is called  an embedded curve \emph{of Mumford type}.
All other embeddings come
from direct sums of such representations. \qed
\endproof

   \begin{rmq} A different proof of this classification can be found in \cite{carshcurves}.
 \end{rmq}
 

\subsection{Relation with the Arakelov bounds}\label{sec:Fin}

Let me summarize the pertinent results so far. First of all, the base curve of  any family of $g$--dimensional
polarized Abelian varieties over  a curve attaining the Arakelov bound must be a Shimura curve embedded in $\AA_g$ as
 a rigid curve of Satake type or of Mumford type.
 If such a curve is also the base curve of a family of  curves of genus $g$, then $g\le 4$.  

Conversely,  an embedded  Shimura curve $C_0\subset \AA_g$    of the above type and which is  generically contained in
Torelli locus need not carry a family of semi-stable curves over it, since $\AA_g$ is not a fine moduli space. Passing to
a suitable fine moduli space $\AA_{g,m}$ destroys the injectivity of the period map: it becomes a branched double cover
which is precisely branched in the intersection of the  locus of the hyperelliptic curves and the Torelli locus.  For curves
not meeting this locus, it of course still follows that $g\le 4$.  However,  one cannot a priori  exclude 
special curves $C_0\subset \AA_{g,m}$ that meet the hyperelliptic locus.

Let me give a brief  \emph{outline the proof of  Theorem (B)}  as stated in the Introduction. 
Let $C$ be the Zariski closure of $C_0$ in the moduli space  $\overline{\MM}_{g,m}$
of semi-stable curves and set 
\[
\HH_g=\text{  the closure of the hyperelliptic locus
 in } \overline{\MM}_{g.m}.
 \]
  Now either $C\subset \HH_g$
or $C$ meets the hyperelliptic locus in a finite set of points. 

For the first case, see \cite{shimura}.
In the second case, form the  double cover $\pi: B\to C$ branched in the points of $\HH_g \cap C$.
 This produces a surface $f: X\to B$ fibered in genus $g$ curves. A little care is needed here: the fiber  over
a  point in the ramification locus might be \emph{singular}, but at most of compact type. As we have
seen, this is because the closure of the  Torelli locus in $\AA_{g.m}$ might contain singular curves of compact type. 

On $B$ the Arakelov bound is then 
 no longer  valid: by \cite[Proposition 1.9]{shimura2},  one has to replace  it by
 \begin{equation*}
 d(f)=  \deg \HH^{1,0}_f = 
\deg f_* \omega_{S/B}= \frac 1 2 \cdot     g \cdot( \deg \Omega^1_B(\log (\Sigma) )- h_0), \quad
h_0= |\HH_g\cap C_0|.
\end{equation*} 
This follows directly from the Hurwitz formula for the double cover $\pi$. One can use this  to rewrite 
 \eqref{eqn:Ineq2} as
\begin{equation*}
\label{eqn:Ineq3}
 4\left(1-\frac 1g\right) +  \frac{1}{d(f)} \cdot  \left( 2(g-1) h_0+
 \sum_{s\in \Sigma_c}   3\delta_s  -\beta(R) -\sum_{j\in J} \beta(f_j)\right)\ge \lambda(f)= \frac{\omega_{X/B}^2}{d(f)}.
\end{equation*}
However, one needs a better inequality. To state it conveniently,   some further notation is used: 
\[
\aligned
n_{j,s} &= \# \text{  of components of genus $j$ in $X_s$,}\\
\epsilon_s &=  \left(\sum_{j\ge 1}   n_{j,s} \right)-1   \\
\theta_s& = \left(\sum_{j\ge 2}3n_{j,s}+2 n_{1,s}\right) -3   \\
\Lambda_0& = \HH_g\cap C_0.
\endaligned
\]
Then  by \cite[Theorem 1.13]{shimura2} one has:
\begin{equation}
\label{eqn:Ineq3}
 4\left(1-\frac 1g\right) +   \frac{2(g-1)h_0}{d(f)}  +\frac{1}{d(f)}\cdot   \left(  
\frac 3 2 \sum_{s\in \Sigma_c\cap \Lambda_0}   \epsilon_s   + \sum_{s\in \Sigma_c\setminus \Lambda_0 }\theta_s
 \right) \ge \lambda(f)= \frac{\omega_{X/B}^2}{d(f)}.
\end{equation}
This inequality turns out to be  a consequence of the refined BMY-inequality by pulling back the family along a   carefully chosen   ramified cover of the base curve $B$ and replacing  $R$ by another curve.

In order to estimate  $h_0$, the idea is to
replace the Moriwaki inequality by a different  slope inequality in which  $h_0$ also   appears.
%
 This  refined slope inequality 
from    \cite[Theorem 1.14]{shimura2}  is a  highly non-trivial
 result which depends on an analysis of a morphism between auxiliary Higgs bundles and a fine analysis
of the singularities of the hyperelliptic singular fibers. It  is valid as soon as $g\ge 7$ and reads
as follows:
 \begin{equation}\label{eqn:Ineq4}
\lambda(f)\ge  \frac{5g-6}{g}  +\frac{1}{d(f)} \cdot \left[ 2(g-2) h_0 \right]   
                   +\frac{1}{d(f)} \cdot \left[\sum_{s\in \Sigma_c\cap \Lambda_0} 
2\epsilon_s  +\sum_{s\in \Sigma_c\setminus \Lambda_0 }\theta_s\right]. 
\end{equation}
Comparing  with \eqref{eqn:Ineq3}, one gets a bound for $h_0$:
\[
\frac{2  h_0}{d(f)} \ge \frac{g-2  }{g} +
\frac{1} {d(f)} \cdot \left[ \sum_{s\in \Sigma_c\cap \Lambda_0}\half \epsilon_s \right]\ge  \frac{g-2  }{g}.
\]
Here one uses that   $\epsilon_s=(\sum_{j\ge 1}  n_{j,s})-1\ge 0$ since there must be non-rational components in the fiber $X_s$. 

Now use  this bound together with the primitive upperbound \eqref{eqn:xiao}, together with \eqref{eqn:Ineq4} as follows:
\[
12 \ge \lambda(f) \ge \frac{5g-6}{g} + \frac{(g-2)^2}{g} +  \text{ remaining term.}
\]
where the remaining  term  is  the  \emph{non-negative} expression
$$
\frac{1}{d(f)} \cdot \left[\sum_{s\in \Sigma_c\cap \Lambda_0} 
2\epsilon_s  +\sum_{s\in \Sigma_c\setminus \Lambda_0 }\theta_s\right].
$$
 It follows that 
\[
0\ge \frac{5g-6}{g} + \frac{(g-2)^2}{g} -12 = \frac{g^2-11g -2}{g} \implies g\le 11.\qed
\]

\begin{rmks}
1) Non-rigidly embedded curves also have been classified; in fact  in \cite{addington,Saitoiii} one finds a
 classification  scheme for all special subvarieties of $\AA_g$. The variation of Hodge structures  corresponding to 
 non-rigid curve embeddings   might be 
 decomposable over $\Q$ and the  corresponding Higgs bundles may contain unitary subsystems. 
 The number of possibilities depend on $g$ and the field $F$ over which the quaternion algebra is defined.
 As  $g$  grows and the index $[F:\Q]$ grows, the number of possibilities grows astronomically.\\
 2) Let $C_0=C\setminus \Delta$ be  a  curve  
contained in the Torelli locus.  By the original Arakelov result, any such curve is rigidly embedded in $\MM_g$
and there are only finitely many of these due to the  boundedness statement in Arakelov's theorem.  So, fixing
$(C,\Delta)$ with $C_0$ a Shimura curve, there are only finitely many embeddings of $C_0$ whose image is
generically contained in the Torelli locus. In particular, among the infinitely many possible embeddings $C_0\into \AA_g$, only finitely
many can be generically contained in the Torelli locus. In particular, if $[F:\Q]$ is large, such a curve is not contained in it.
\\
 3) Jacobians of smooth curves have an irreducible polarization
which give  $\Q$--irreducible variations of Hodge structures. However, the corresponding Higgs bundles (which only "see" the $\R$-structure) can decompose and hence could in particular contain unitary Higgs subbundles defined over some
real number field $\not=\Q$. Note that for such a 
subbundle the monodromy need no longer be finite.
Indeed, this can be seen to happen for the examples in \cite{Torloc}. See also \cite{vbvhs}.
\end{rmks}

\end{document}